\documentclass[reqno,onecolumn,12pt]{amsart}
\usepackage{amsfonts}
\usepackage{amssymb}
\usepackage{amsmath}
\usepackage{graphicx}
\usepackage{amscd}

\newtheorem{theorem}{Theorem}
\theoremstyle{plain}

\newtheorem{definition}{Definition}

\numberwithin{equation}{section}

\input{tcilatex}

\begin{document}
\author{}
\title{}
\maketitle

\begin{center}
\thispagestyle{empty} \textbf{A\ NEW\ APPROACH\ TO\ MODIFIED }$q$\textbf{%
-BERNSTEIN\ POLYNOMIALS FOR\ FUNCTIONS OF TWO VARIABLES WITH THEIR
GENERATING\ AND INTERPOLATION\ FUNCTIONS}

\bigskip

\textbf{Mehmet ACIKGOZ }and\textbf{\ Serkan ARACI}

University of Gaziantep, Faculty of Arts and Science, Department of
Mathematics, 27310 Gaziantep, Turkey

acikgoz@gantep.edu.tr; mtsrkn@hotmail.com\bigskip

\textbf{{\Large {Abstract}}}
\end{center}

The aim of this paper is to give a new approach to modified $q$-Bernstein
polynomials for functions of two variables. By using these type polynomials,
we derive recurrence formulas and some new interesting identities related to
the second kind Stirling numbers and generalized Bernoulli polynomials.
Moreover, we give the generating function and interpolation function of
these modified $q$-Bernstein polynomials of two variables and also give the
derivatives of these polynomials and their generating function.

\bigskip



\section{Introduction, Definitions and Notations}

In approximation theory, the Bernstein polynomials, named after their
creater S. N. Bernstein in 1912, have been studied by many researchers for a
long time. But nothing about generating function of Bernstein polynomials
were known in the literature. Recently, Simsek and Acikgoz, (\cite%
{SimsekAcikgoz1}), constructed a new generating function of ($q$-) Bernstein
type polynomials based on the $q$-analysis. They gave some new relations
related to ($q$-) Bernstein type polynomials, Hermite polynomials, Bernoulli
polynomials of higher order and the second kind Stirling numbers. By
applying Mellin transformation to this generating function they defined the
interpolation function of ($q$-) Bernstein type polynomials. They gave some
relations and identities on these polynomials. They constructed the
generating function for classical Bernstein polynomials, and for Bernstein
polynomials for functions of two variables and gave their properties (see 
\cite{Acikgoz1}, \cite{Acikgoz2}, \cite{Acikgoz3}, for details).

Throughout this paper, we use some notations like $%
\mathbb{N}
,$ $%
\mathbb{N}
_{0}$ and $D,$ where $%
\mathbb{N}
$ denotes the set of natural numbers, $%
\mathbb{N}
_{0}:=%
\mathbb{N}
\tbigcup \left\{ 0\right\} $ and $D=\left[ 0,1\right] $.

Let $C\left( D\times D\right) $ denotes the set of continuous functions on $%
D $. For $f\in C\left( D\times D\right) $%
\begin{eqnarray}
\mathbf{B}_{n,m}\left( f;x,y\right)
&:&=\dsum\limits_{k=0}^{n}\dsum\limits_{j=0}^{m}f\left( \frac{k}{n},\frac{j}{%
m}\right) \binom{n}{k}\binom{m}{j}x^{k}y^{j}\left( 1-x\right) ^{n-k}\left(
1-y\right) ^{m-j}  \notag \\
&=&\dsum\limits_{k=0}^{n}\dsum\limits_{j=0}^{m}f\left( \frac{k}{n},\frac{j}{m%
}\right) B_{k,j;n,m}\left( x,y\right)  \label{Equation1}
\end{eqnarray}%
where $\binom{n}{k}=\frac{n\left( n-1\right) \cdots \left( n-k+1\right) }{k!}%
.$ Here $\mathbf{B}_{n,m}\left( f;x,y\right) $ is called the Bernstein
operator of two variables of order $n+m$ for $f$. For $k,j,n,m\in 
\mathbb{N}
_{0}$, the Bernstein polynomial of two variables of degree $n+m$ is defined
by%
\begin{equation}
B_{k,j;n,m}\left( x,y\right) =\binom{n}{k}\binom{m}{j}x^{k}y^{j}\left(
1-x\right) ^{n-k}\left( 1-y\right) ^{m-j},  \label{Equation2}
\end{equation}%
where $x\in D$ and $y\in D$. Thus, throughout this work, we will assume that 
$x\in D$ and $y\in D$. Then, we easily see the following%
\begin{equation}
B_{k,j;n,m}\left( x,y\right) =B_{k,n}\left( x,y\right) B_{j,m}\left(
x,y\right)  \label{Equation3}
\end{equation}%
and they form a partition of unity; that is; 
\begin{equation}
\dsum\limits_{k=0}^{n}\dsum\limits_{j=0}^{m}B_{k,j;n,m}\left( x,y\right) =1
\label{Equation4}
\end{equation}%
and by using the definition of Bernstein polynomials for functions of two
variables, it is not difficult to prove the property given above as%
\begin{equation}
\dsum\limits_{k=0}^{n}\dsum\limits_{j=0}^{m}B_{k,n}\left( x,y\right)
B_{j,m}\left( x,y\right) =1.  \label{Equation5}
\end{equation}%
Some Bernstein polynomials of two variables are given below:%
\begin{equation*}
B_{0,0;1,0}\left( x,y\right) =\left( 1-x\right) ,\text{ }B_{0,0;0,1}\left(
x,y\right) =\left( 1-y\right) ,\text{ }B_{0,0;1,1}\left( x,y\right) =\left(
1-x\right) \left( 1-y\right) ,
\end{equation*}%
\begin{equation*}
B_{0,1;1,1}\left( x,y\right) =y\left( 1-x\right) ,\text{ }B_{1,0;1,1}\left(
x,y\right) =x\left( 1-y\right) ,\text{ }B_{1,1;1,1}\left( x,y\right) =xy%
\text{.}
\end{equation*}%
Also, $B_{k,j;n,m}\left( x,y\right) =0$ for $k>n$ or $j>m$, because $\binom{n%
}{k}=0$ or $\binom{m}{j}=0.~$There are $nm+n+m+1,\ n+m$-th degree Bernstein
polynomials (see \cite{Acikgoz3} and \cite{Buyukyazici2} for details).

Some researchers have used the Bernstein polynomials of two variables in
approximation theory (See \cite{Buyukyazici1}, \cite{Buyukyazici2}). But no
result was known anything about the generating function of these
polynomials. Note that for $k,j,n,m\in 
\mathbb{N}
_{0}$, we have%
\begin{eqnarray*}
\frac{\left( tx\right) ^{k}\left( ty\right) ^{j}e^{2t}}{k!j!e^{t\left(
x+y\right) }} &=&\frac{t^{k}x^{k}t^{j}y^{j}}{k!j!}e^{t\left( 1-x\right)
}e^{t\left( 1-y\right) } \\
&=&\frac{x^{k}}{k!}\left( t^{k}\dsum\limits_{n=0}^{\infty }\frac{\left(
1-x\right) ^{n}}{n!}t^{n}\right) \frac{y^{j}}{j!}\left(
t^{j}\dsum\limits_{m=0}^{\infty }\frac{\left( 1-y\right) ^{m}}{m!}%
t^{m}\right) \\
&=&\dsum\limits_{n=k}^{\infty }\dsum\limits_{m=j}^{\infty }B_{k,j;n,m}\left(
x,y\right) \frac{t^{n}}{n!}\frac{t^{m}}{m!}
\end{eqnarray*}%
From the above, we obtain the generating function for $B_{k,j;n,m}\left(
x,y\right) $ as follows:%
\begin{equation}
F_{k,j}\left( t;x,y\right) =\frac{\left( tx\right) ^{k}\left( ty\right)
^{j}e^{t\left( 2-\left( x+y\right) \right) }}{k!j!}=\dsum\limits_{n=k}^{%
\infty }\dsum\limits_{m=j}^{\infty }B_{k,j;n,m}\left( x,y\right) \frac{t^{n}%
}{n!}\frac{t^{m}}{m!},  \label{Equation55}
\end{equation}%
where $k,j,n,m\in 
\mathbb{N}
_{0}$. We notice that,

\begin{equation*}
B_{k,j;n,m}\left( x,y\right) =\left\{ 
\begin{array}{cccccc}
\binom{n}{k}\binom{m}{j}x^{k}y^{j}\left( 1-x\right) ^{n-k}\left( 1-y\right)
^{m-j} & , & \text{if} & n\geq k & \text{and} & m\geq j \\ 
0 & , & \text{if} & n<k & \text{or} & m<j%
\end{array}%
\right.
\end{equation*}%
for $n,k,m,j\in 
\mathbb{N}
_{0}$ (for details, see \cite{Acikgoz2}).

Let $q\in \left( 0,1\right) $. Then, $q$-integer of $x$ by $[x]_{q}:=\frac{%
1-q^{x}}{1-q}$ \ and \ $[x]_{-q}:=\frac{1-\left( -q\right) ^{x}}{1+q}$ ( See 
\cite{Kim-Lee-Chae Jang}-\cite{SimsekAcikgoz1} for details). Note that $%
\underset{q\rightarrow 1}{\lim }[x]_{q}=x$. \cite{Kim-Lee-Chae Jang}
motivated the authors to write this paper and we have extended the results
given in that paper to modified $q$-Bernstein polynomials of two variables.

\section{The Modified $q$-Bernstein Polynomials for Functions of two
Variables}

For $0\leq k\leq n$ and $0\leq j\leq m$, the $q$-Bernstein polynomials of
degree $n+m$ are defined by%
\begin{equation}
B_{k,j;n,m}\left( x,y;q\right) =\left\{ 
\begin{array}{cccccc}
\binom{n}{k}\binom{m}{j}[x]_{q}^{k}[y]_{q}^{j}[1-x]_{q}^{n-k}[1-y]_{q}^{m-j}
& , & \text{if} & n\geq k & \text{and} & m\geq j \\ 
0 & , & \text{if} & n<k & \text{or} & m<j%
\end{array}%
\right. .  \label{Equation56}
\end{equation}

For $q\in \left( 0,1\right) $, consider the $q$-extension of \ (\ref%
{Equation55}) as follows:%
\begin{eqnarray}
F_{k,j}\left( t,q;x,y\right) &=&\frac{\left( t[x]_{q}\right) ^{k}\left(
t[y]_{q}\right) ^{j}}{k!j!}e^{t\left( [1-x]_{q}+[1-y]_{q}\right) }  \notag \\
&=&\frac{[x]_{q}^{k}}{k!}\left( \dsum\limits_{n=0}^{\infty }\frac{%
[1-x]_{q}^{n}}{n!}t^{n+k}\right) \frac{[y]_{q}^{j}}{j!}\left(
\dsum\limits_{m=0}^{\infty }\frac{[1-y]_{q}^{m}}{m!}t^{m+j}\right)  \notag \\
&=&\dsum\limits_{n=k}^{\infty }\dsum\limits_{m=j}^{\infty }\binom{n}{k}%
\binom{m}{j}[x]_{q}^{k}[y]_{q}^{j}[1-x]_{q}^{n-k}[1-y]_{q}^{m-j}\frac{t^{n}}{%
n!}\frac{t^{m}}{m!}  \label{Equation900}
\end{eqnarray}%
where $k,j,n,m\in 
\mathbb{N}
_{0}$. Note that $\underset{q\rightarrow 1}{\lim }F_{k,j}\left(
t,q:x,y\right) =F_{k,j}\left( t;x,y\right) .$

\begin{definition}
The modified $q$-Bernstein polynomials for functions of two variables is
defined by means of the following generating function: 
\begin{equation}
F_{k,j}\left( t,q;x,y\right) =\frac{\left( t[x]_{q}\right) ^{k}\left(
t[y]_{q}\right) ^{j}}{k!j!}e^{t\left( [1-x]_{q}+[1-y]_{q}\right)
}=\dsum\limits_{n=0}^{\infty }\dsum\limits_{m=0}^{\infty }B_{k,j;n,m}\left(
x,y;q\right) \frac{t^{n}}{n!}\frac{t^{m}}{m!}  \label{Equation901}
\end{equation}%
where $k,j,n,m\in 
\mathbb{N}
_{0}$.
\end{definition}

By comparing the coefficients of (\ref{Equation900}) and (\ref{Equation901}%
), we obtain a formula for modified $q$-Bernstein polynomials of two
variables given in the following theorem:

\begin{theorem}
For $k,j,n,m\in 
\mathbb{N}
_{0}$, then, we have%
\begin{equation}
B_{k,j;n,m}\left( x,y;q\right) =\left\{ 
\begin{array}{cccccc}
\binom{n}{k}\binom{m}{j}[x]_{q}^{k}[y]_{q}^{j}[1-x]_{q}^{n-k}[1-y]_{q}^{m-j}
& , & \text{if} & n\geq k & \text{and} & m\geq j \\ 
0 & , & \text{if} & n<k & \text{or} & m<j%
\end{array}%
\right. .  \label{Equation59}
\end{equation}
\end{theorem}

\begin{theorem}
\textbf{(Recurrence Formula for }$B_{k,j;n,m}\left( x,y;q\right) $\textbf{)}
For $k,j,n,m\in 
\mathbb{N}
_{0}$, we have%
\begin{eqnarray*}
B_{k,j;n,m}\left( x,y;q\right) &=&[1-x]_{q}[1-y]_{q}B_{k,j;n-1,m-1}\left(
x,y;q\right) +[1-x]_{q}[y]_{q}B_{k,j-1;n-1,m-1}\left( x,y;q\right) \\
&&+[x]_{q}[1-y]_{q}B_{k-1,j;n-1,m-1}\left( x,y;q\right)
+[x]_{q}[y]_{q}B_{k-1,j-1;n-1,m-1}\left( x,y;q\right) .
\end{eqnarray*}

\begin{proof}
By using the definition of Bernstein polynomials for functions of two
variables defined by (\ref{Equation59}), we have 
\begin{eqnarray*}
B_{k,j;n,m}\left( x,y;q\right) &=&\binom{n}{k}\binom{m}{j}%
[x]_{q}^{k}[y]_{q}^{j}[1-x]_{q}^{n-k}[1-y]_{q}^{m-j} \\
&=&\left[ \binom{n-1}{k}+\binom{n-1}{k-1}\right] \left[ \binom{m-1}{j}+%
\binom{m-1}{j-1}\right] [x]_{q}^{k}[y]_{q}^{j}[1-x]_{q}^{n-k}[1-y]_{q}^{m-j}
\\
&=&[1-x]_{q}[1-y]_{q}B_{k,j;n-1,m-1}\left( x,y;q\right)
+[1-x]_{q}[y]_{q}B_{k,j-1;n-1,m-1}\left( x,y;q\right) \\
&&+[x]_{q}[1-y]_{q}B_{k-1,j;n-1,m-1}\left( x,y;q\right)
+[x]_{q}[y]_{q}B_{k-1,j-1;n-1,m-1}\left( x,y;q\right) .
\end{eqnarray*}
\end{proof}
\end{theorem}

\begin{theorem}
For $k,j,n,m\in 
\mathbb{N}
_{0}$, we get 
\begin{equation}
B_{n-k,m-j;n,m}\left( 1-x,1-y;q\right) =B_{k,j;n,m}\left( x,y;q\right)
\label{Equation151}
\end{equation}%
and%
\begin{equation*}
\mathbf{B}_{n,m}\left( 1:x,y,q\right) =\left( 1+\left( 1-q\right)
[x]_{q}[1-x]_{q}\right) ^{n}\times \left( 1+\left( 1-q\right)
[y]_{q}[1-y]_{q}\right) ^{m}.
\end{equation*}

\begin{proof}
Let $f$ be a continuous function of two variables on $D\times D$. Then the
modified $q$-Bernstein operator of order $n+m$ for $f$ is defined by 
\begin{equation}
\mathbf{B}_{n,m}\left( f:x,y,q\right)
=\dsum\limits_{k=0}^{n}\dsum\limits_{j=0}^{m}f\left( \frac{k}{n},\frac{j}{m}%
\right) B_{k,j;n,m}\left( x,y;q\right)  \label{Equation160}
\end{equation}%
where $0\leq x\leq 1,$ $0\leq y\leq 1$, $n,m\in 
\mathbb{N}
.$ From Theorem 1 and the definition of modified $q$-Bernstein operator
given by (\ref{Equation160}) for $f\left( x,y\right) =xy$, we have 
\begin{eqnarray*}
\mathbf{B}_{n,m}\left( f:x,y,q\right)
&=&\dsum\limits_{k=0}^{n}\dsum\limits_{j=0}^{m}f\left( \frac{k}{n},\frac{j}{m%
}\right) \binom{n}{k}\binom{m}{j}%
[x]_{q}^{k}[1-x]_{q}^{n-k}[y]_{q}^{k}[1-y]_{q}^{m-j} \\
&=&[x]_{q}\left( 1-[1-x]_{q}[x]_{q}\left( q-1\right) \right) ^{n-1}\times
\lbrack y]_{q}\left( 1-[1-y]_{q}[y]_{q}\left( q-1\right) \right) ^{m-1} \\
&=&f\left( [x]_{q},[y]_{q}\right) \left( 1+\left( 1-q\right)
[x]_{q}[1-x]_{q}\right) ^{n-1}\times \left( 1+\left( 1-q\right)
[y]_{q}[1-y]_{q}\right) ^{m-1}
\end{eqnarray*}%
From Theorem 1, we have 
\begin{eqnarray*}
\mathbf{B}_{n,m}\left( 1:x,y,q\right)
&=&\dsum\limits_{k=0}^{n}\dsum\limits_{j=0}^{m}B_{k,j;n,m}\left( x,y;q\right)
\\
&=&\dsum\limits_{k=0}^{n}B_{k,n}\left( x,q\right)
\dsum\limits_{j=0}^{m}B_{j,m}\left( y,q\right) \\
&=&\left( 1+\left( 1-q\right) [x]_{q}[1-x]_{q}\right) ^{n}\left( 1+\left(
1-q\right) [y]_{q}[1-y]_{q}\right) ^{m}.
\end{eqnarray*}%
The modified $q$-Bernstein polynomials of two variables are symmetric
polynomials:%
\begin{eqnarray*}
B_{n-k,m-j;n,m}\left( 1-x,1-y;q\right) &=&\binom{n}{n-k}%
[x]_{q}^{k}[1-x]_{q}^{n-k}\binom{n}{m-j}[y]_{q}^{j}[1-y]_{q}^{m-j} \\
&=&\binom{n}{k}[1-x]_{q}^{k}[x]_{q}^{n-k}\binom{n}{j}%
[1-y]_{q}^{j}[y]_{q}^{m-j} \\
&=&B_{k,j;n,m}\left( x,y;q\right) .
\end{eqnarray*}%
by replacing $k$ by $n-k$ and $j$ by $m-j$.
\end{proof}
\end{theorem}

\begin{theorem}
For $\xi ,\rho \in 
\mathbb{C}
$, and for $n,m\in 
\mathbb{N}
$, then, we procure 
\begin{equation}
B_{k,j;n,m}\left( x,y;q\right) =-\frac{n!m!}{4\pi ^{2}}\doint\limits_{C}%
\doint\limits_{C}\frac{\left( \left[ x\right] _{q}\xi \right) ^{k}\left( %
\left[ y\right] _{q}\rho \right) ^{j}}{k!j!}e^{\left( \left[ 1-x\right]
_{q}\xi +\left[ 1-y\right] _{q}\rho \right) }\frac{d\xi }{\xi ^{n+1}}\frac{%
d\rho }{\rho ^{m+1}}  \label{Equation153}
\end{equation}%
where $C$ is a circle around the origin and integration is in the positive
direction.

\begin{proof}
By using the definition of the modified $q$-Bernstein polynomials of two
variables and the basic theory of complex analysis including Laurent series
that%
\begin{eqnarray}
&&\doint\limits_{C}\doint\limits_{C}\frac{\left( \left[ x\right] _{q}\xi
\right) ^{k}\left( \left[ y\right] _{q}\rho \right) ^{j}}{k!j!}e^{\left( %
\left[ 1-x\right] _{q}\xi +\left[ 1-y\right] _{q}\rho \right) }\frac{d\xi }{%
\xi ^{n+1}}\frac{d\rho }{\rho ^{m+1}}  \notag \\
&=&\dsum\limits_{l=0}^{\infty }\dsum\limits_{r=0}^{\infty
}\doint\limits_{C}\doint\limits_{C}\frac{B_{k,l}\left( x,q\right) \xi ^{l}}{%
l!}\frac{B_{j,r}\left( y,q\right) \rho ^{r}}{r!}\frac{d\xi }{\xi ^{n+1}}%
\frac{d\rho }{\rho ^{m+1}}  \label{Equation152} \\
&=&\left( 2\pi i\right) ^{2}\left( \frac{B_{k,j;n,m}\left( x,y;q\right) }{%
n!m!}\right) \text{.}  \notag
\end{eqnarray}%
By using (\ref{Equation153} ) and (\ref{Equation152}), we obtain 
\begin{equation*}
\frac{n!m!}{\left( 2\pi i\right) ^{2}}\doint\limits_{C}\doint\limits_{C}%
\frac{\left( \left[ x\right] _{q}\xi \right) ^{k}}{k!}\frac{\left( \left[ y%
\right] _{q}\rho \right) ^{j}}{j!}\frac{d\xi }{\xi ^{n+1}}\frac{d\rho }{\rho
^{m+1}}=B_{k,j;nm}\left( x,y;q\right)
\end{equation*}%
and 
\begin{eqnarray}
&&\doint\limits_{C}\doint\limits_{C}\frac{\left( \left[ x\right] _{q}\xi
\right) ^{k}}{k!}\frac{\left( \left[ y\right] _{q}\rho \right) ^{j}}{j!}%
e^{\left( \left[ 1-x\right] _{q}\xi +\left[ 1-y\right] _{q}\rho \right) }%
\frac{d\xi }{\xi ^{n+1}}\frac{d\rho }{\rho ^{m+1}}  \label{Equation154} \\
&=&\left( 2\pi i\right) ^{2}\left( \frac{\left[ x\right] _{q}^{k}\left[ y%
\right] _{q}^{j}\left[ 1-x\right] _{q}^{n-k}\left[ 1-y\right] _{q}^{m-j}}{%
k!j!\left( n-k\right) !\left( m-j\right) !}\right) \text{.}  \notag
\end{eqnarray}%
We also obtain from (\ref{Equation151}) and (\ref{Equation154}) that 
\begin{eqnarray}
&&\frac{n!m!}{\left( 2\pi i\right) ^{2}}\doint\limits_{C}\doint\limits_{C}%
\frac{\left( \left[ x\right] _{q}\xi \right) ^{k}\left( \left[ y\right]
_{q}\rho \right) ^{j}}{k!j!}e^{\left( \left[ 1-x\right] _{q}\xi +\left[ 1-y%
\right] _{q}\rho \right) }\frac{d\xi }{\xi ^{n+1}}\frac{d\rho }{\rho ^{m+1}}
\label{Equation155} \\
&=&\binom{n}{k}\binom{m}{j}\left[ x\right] _{q}^{k}\left[ 1-x\right]
_{q}^{n-k}\left[ y\right] _{q}^{j}\left[ 1-y\right] _{q}^{m-j}.  \notag
\end{eqnarray}%
Therefore we see that from (\ref{Equation152}) and (\ref{Equation155}) that 
\begin{equation*}
B_{k,j;n,m}\left( x,y;q\right) =\binom{n}{k}\binom{m}{j}\left[ x\right]
_{q}^{k}\left[ 1-x\right] _{q}^{n-k}\left[ y\right] _{q}^{j}\left[ 1-y\right]
_{q}^{m-j}.
\end{equation*}
\end{proof}
\end{theorem}

\begin{theorem}
(\textbf{The Derivative Formula for }$B_{k,j;n,m}\left( x,y;q\right) $) For $%
k,j,n,m\in 
\mathbb{N}
$, then, we derive the following%
\begin{eqnarray*}
\frac{\partial ^{2}}{\partial x\partial y}\left( B_{k,j;n,m}\left(
x,y;q\right) \right) &=&nm(q^{x+y}B_{k-1,j-1;n-1,m-1}\left( x,y;q\right)
-q^{x-y+1}B_{k-1,j;n-1,m-1}\left( x,y;q\right) \\
&&-q^{1-x+y}B_{k,j-1;n-1,m-1}\left( x,y;q\right) +q^{2-\left( x+y\right)
}B_{k,j;n-1,m-1}\left( x,y;q\right) )\frac{\ln ^{2}q}{\left( q-1\right) ^{2}}%
.
\end{eqnarray*}

\begin{proof}
Using the definition of modified $q$-Bernstein polynomials for functions of
two variables and the property (\ref{Equation3}), we have 
\begin{equation*}
\frac{\partial ^{2}}{\partial x\partial y}\left( B_{k,j;n,m}\left(
x,y;q\right) \right) =\frac{\partial ^{2}}{\partial x\partial y}\left(
B_{k,n}\left( x;q\right) B_{j,m}\left( y;q\right) \right) =\frac{d}{dx}%
\left( B_{k,n}\left( x;q\right) \right) \frac{d}{dy}\left( B_{j,m}\left(
y;q\right) \right)
\end{equation*}%
and after some calculations, the proof is complete.
\end{proof}
\end{theorem}

Therefore, we can write the modified $q$-Bernstein polynomials for functions
of two variables as a linear combination of polynomials of higher order as
follows:

\begin{theorem}
For $k,j,n,m\in 
\mathbb{N}
_{0}$, we have%
\begin{eqnarray*}
&&\left( 1+\left( 1-q\right) [x]_{q}\left[ 1-x\right] _{q}\right) \left(
1+\left( 1-q\right) [y]_{q}\left[ 1-y\right] _{q}\right) B_{k,j;n,m}\left(
x,y;q\right) \\
&=&\left( \frac{n-k+1}{n+1}\right) \left( \frac{m-j+1}{m+1}\right)
B_{k,j;n+1,m+1}\left( x,y;q\right) +\left( \frac{n-k+1}{n+1}\right) \left( 
\frac{j+1}{m+1}\right) B_{k,j+1;n+1,m+1}\left( x,y;q\right) \\
&&+\left( \frac{k+1}{n+1}\right) \left( \frac{m-j+1}{m+1}\right)
B_{k+1,j;n+1,m+1}\left( x,y;q\right) +\left( \frac{k+1}{n+1}\right) \left( 
\frac{j+1}{m+1}\right) B_{k+1,j+1;n+1,m+1}\left( x,y;q\right) .
\end{eqnarray*}

\begin{proof}
It follows after expanding the series and some algebraic operations.
\end{proof}
\end{theorem}

\begin{theorem}
For $k,j,n,m\in 
\mathbb{N}
_{0}$, we have%
\begin{equation*}
B_{k,j;n,m}\left( x,y;q\right) =\left( \frac{n-k+1}{k}\right) \left( \frac{%
m-j+1}{j}\right) \left( \frac{[x]_{q}[y]_{q}}{[1-x]_{q}[1-y]_{q}}\right)
B_{k-1,j-1;n,m}\left( x,y;q\right) .
\end{equation*}

\begin{proof}
To prove this theorem, we start with the right hand side:%
\begin{eqnarray*}
&&\left( \frac{n-k+1}{k}\right) \left( \frac{m-j+1}{j}\right) \left( \frac{%
[x]_{q}[y]_{q}}{[1-x]_{q}[1-y]_{q}}\right) B_{k-1,j-1;n,m}\left( x,y;q\right)
\\
&=&\frac{n!}{\left( n-k\right) !k!}.\frac{m!}{\left( m-j\right) !j!}\left( 
\frac{[x]_{q}[y]_{q}}{[1-x]_{q}[1-y]_{q}}\right)
[x]_{q}^{k-1}[y]_{q}^{j-1}[1-x]_{q}^{n-k+1}[1-y]_{q}^{m-j+1} \\
&=&\binom{n}{k}\binom{m}{j}%
[x]_{q}^{k}[y]_{q}^{j}[1-x]_{q}^{n-k}[1-y]_{q}^{m-j}=B_{k,j;n,m}\left(
x,y;q\right) .
\end{eqnarray*}
\end{proof}
\end{theorem}

\begin{theorem}
For $k,j,n,m\in 
\mathbb{N}
_{0}$, we obtain%
\begin{equation*}
B_{k,j;n,m}\left( x,y;q\right) =\dsum\limits_{l=k}^{n}\dsum\limits_{r=j}^{m}%
\binom{n}{l}\binom{l}{k}\binom{m}{j}\binom{j}{r}\left( -1\right)
^{l-k+r-j}q^{\left( l-k\right) \left( 1-x\right) +\left( r-j\right) \left(
1-y\right) }[x]_{q}^{l}[y]_{q}^{r}.
\end{equation*}

\begin{proof}
From the definition of modified $q$-Bernstein polynomials of two variables
and binomial theorem with $k,j,n,m\in 
\mathbb{N}
_{0}$, we have%
\begin{eqnarray*}
B_{k,j;n,m}\left( x,y;q\right) &=&\binom{n}{k}\binom{m}{j}%
[x]_{q}^{k}[1-x]_{q}^{n-k}[y]_{q}^{j}[1-y]_{q}^{m-j} \\
&=&\binom{n}{k}\binom{m}{j}[x]_{q}^{k}[y]_{q}^{j}\left( 1-q^{1-x}[x]\right)
^{n-k}\left( 1-q^{1-y}[y]\right) ^{m-j} \\
&=&\dsum\limits_{l=k}^{n}\dsum\limits_{r=j}^{m}\binom{n}{l}\binom{l}{k}%
\binom{m}{j}\binom{j}{r}\left( -1\right) ^{l-k+r-j}q^{\left( l-k\right)
\left( 1-x\right) +\left( r-j\right) \left( 1-y\right)
}[x]_{q}^{l}[y]_{q}^{r}.
\end{eqnarray*}
\end{proof}
\end{theorem}

\begin{theorem}
The following identity 
\begin{equation*}
\left( \lbrack x]_{q}[y]_{q}\right) ^{l}=\frac{1}{\left( \left[ 1-x\right]
_{q}+\left[ x\right] _{q}\right) ^{n-l}\left( \left[ 1-y\right] _{q}+\left[ y%
\right] _{q}\right) ^{m-l}}\dsum\limits_{k=l}^{n}\dsum\limits_{j=l}^{m}\frac{%
\binom{k}{l}\binom{j}{l}}{\binom{n}{l}\binom{m}{l}}B_{k,j;n,m}\left(
x,y;q\right)
\end{equation*}%
is true.

\begin{proof}
We easily see that from the property of the modified $q$-Bernstein
polynomials of two variables that%
\begin{eqnarray*}
\dsum\limits_{k=1}^{n}\dsum\limits_{j=1}^{m}\frac{kj}{nm}B_{k,j;n,m}\left(
x,y;q\right) &=&\dsum\limits_{k=1}^{n}\dsum\limits_{j=1}^{m}\binom{n-1}{k-1}%
\binom{m-1}{j-1}\left[ x\right] _{q}^{k}\left[ y\right] _{q}^{j}\left[ 1-x%
\right] _{q}^{n-k}\left[ 1-y\right] _{q}^{m-j} \\
&=&\left[ x\right] _{q}\left[ y\right] _{q}\left( \left[ x\right] _{q}+\left[
1-x\right] _{q}\right) ^{n-1}\left( \left[ y\right] _{q}+\left[ 1-y\right]
_{q}\right) ^{m-1}
\end{eqnarray*}%
and that 
\begin{eqnarray*}
\dsum\limits_{k=2}^{n}\dsum\limits_{j=2}^{m}\frac{\binom{k}{2}\binom{j}{2}}{%
\binom{n}{2}\binom{m}{2}}B_{k,j;n,m}\left( x,y;q\right)
&=&\dsum\limits_{k=2}^{n}\dsum\limits_{j=2}^{m}\binom{n-2}{k-2}\binom{m-2}{%
j-2}\left[ x\right] _{q}^{k}\left[ y\right] _{q}^{j}\left[ 1-x\right]
_{q}^{n-k}\left[ 1-y\right] _{q}^{m-j} \\
&=&\left[ x\right] _{q}^{2}\left[ y\right] _{q}^{2}\left( \left[ x\right]
_{q}+\left[ 1-x\right] _{q}\right) ^{n-2}\left( \left[ y\right] _{q}+\left[
1-y\right] _{q}\right) ^{m-2}
\end{eqnarray*}%
Continuing this way, we have 
\begin{equation*}
\dsum\limits_{k=l}^{n}\dsum\limits_{j=l}^{m}\frac{\binom{k}{l}\binom{j}{l}}{%
\binom{n}{l}\binom{m}{l}}B_{k,j;n,m}\left( x,y;q\right) =\left[ x\right]
_{q}^{l}\left[ y\right] _{q}^{l}\left( \left[ x\right] _{q}+\left[ 1-x\right]
_{q}\right) ^{n-l}\left( \left[ y\right] _{q}+\left[ 1-y\right] _{q}\right)
^{m-l}
\end{equation*}%
and after some algebraic operations, we obtain the desired result.
\end{proof}
\end{theorem}

We see that from the theorem above, it is possible to write $\left(
[x]_{q}[y]_{q}\right) ^{k}$ as a linear combination of the two variables
modified $q$-Bernstein polynomials.

For $k\in 
\mathbb{N}
_{0}$, the Bernoulli polynomials of degree $k$ are defined by 
\begin{equation*}
\left( \frac{t}{e^{t}-1}\right) ^{k}e^{xt}=\underset{k-times}{\underbrace{%
\left( \frac{t}{e^{t}-1}\right) \times \cdots \times \left( \frac{t}{e^{t}-1}%
\right) }}e^{xt}=\dsum\limits_{n=0}^{\infty }B_{n}^{(k)}\left( x\right) 
\frac{t^{n}}{n!},
\end{equation*}%
and $B_{n}^{\left( k\right) }=B_{n}^{(k)}\left( 0\right) $ are called the $n$%
-th Bernoulli numbers of order $k$. It is well known that the second kind
Stirling numbers are defined by%
\begin{equation}
\frac{\left( e^{t}-1\right) ^{k}}{k!}:=\dsum\limits_{n=0}^{\infty }S\left(
n,k\right) \frac{t^{n}}{n!}  \label{Equation80}
\end{equation}%
for $k\in 
\mathbb{N}
$ (see \cite{Kim-Lee-Chae Jang}). By using the above relations we can give
the following theorem:

\begin{theorem}
For $k,j,n,m\in 
\mathbb{N}
_{0}$, we have%
\begin{eqnarray*}
B_{k,j;n,m}\left( x,y;q\right)
&=&[x]_{q}^{k}[y]_{q}^{j}\dsum\limits_{l=0}^{n}\dsum\limits_{r=0}^{m}\binom{n%
}{l}\binom{m}{r} \\
&&\times B_{l}^{\left( k\right) }\left( \left[ 1-x\right] _{q}\right)
B_{r}^{\left( j\right) }\left( \left[ 1-y\right] _{q}\right) S\left(
n-l,k\right) S\left( m-r,j\right) .
\end{eqnarray*}

\begin{proof}
By using the generating function of modified $q$-Bernstein polynomials of
two variables, we have 
\begin{eqnarray*}
&&\frac{\left( t[x]_{q}\right) ^{k}\left( t[y]_{q}\right) ^{j}}{k!j!}%
e^{t\left( [1-x]_{q}+[1-y]_{q}\right) }=[x]_{q}^{k}[y]_{q}^{j}\left(
\dsum\limits_{n=0}^{\infty }S\left( n,k\right) \frac{t^{n}}{n!}\right)
\left( \dsum\limits_{m=0}^{\infty }S\left( m,j\right) \frac{t^{m}}{m!}\right)
\\
&&\times \left( \dsum\limits_{l=0}^{\infty }B_{l}^{\left( k\right) }\left(
[1-x]_{q}\right) \frac{t^{l}}{l!}\right) \left( \dsum\limits_{r=0}^{\infty
}B_{r}^{\left( j\right) }\left( [1-y]_{q}\right) \frac{t^{r}}{r!}\right) \\
&=&\dsum\limits_{n\geq k}\dsum\limits_{m\geq j}B_{k,j;n,m}\left(
x,y;q\right) \frac{t^{n}}{n!}\frac{t^{m}}{m!}
\end{eqnarray*}
by using the Cauchy product. By comparing last two relations, we have the
desired result.
\end{proof}
\end{theorem}

Let $\Delta $ be the shift difference operator defined by $\Delta f\left(
x\right) =f\left( x+1\right) -f\left( x\right) $. By using the iterative
method we have%
\begin{equation}
\Delta ^{n}f\left( 0\right) =\dsum\limits_{k=0}^{n}\binom{n}{k}\left(
-1\right) ^{n-k}f\left( k\right) ,  \label{Equation90}
\end{equation}%
for $n\in 
\mathbb{N}
$.%
\begin{equation}
\dsum\limits_{n=0}^{\infty }S\left( n,k\right) \frac{t^{n}}{n!}=\frac{1}{k!}%
\dsum\limits_{l=0}^{k}\binom{k}{l}\left( -1\right)
^{k-l}e^{lt}=\dsum\limits_{n=0}^{\infty }\left\{ \frac{1}{k!}%
\dsum\limits_{l=0}^{k}\binom{k}{l}\left( -1\right) ^{k-l}l^{n}\right\} \frac{%
t^{n}}{n!}=\dsum\limits_{n=0}^{\infty }\frac{\Delta ^{k}0^{n}}{k!}\frac{t^{n}%
}{n!}.  \notag
\end{equation}%
By comparing the coefficients on both sides above, we have 
\begin{equation}
S\left( n,k\right) =\frac{\Delta ^{k}0^{n}}{k!}  \label{Equation91}
\end{equation}%
for $n,k\in 
\mathbb{N}
_{0}$. By using the equations (\ref{Equation80}) and (\ref{Equation90}), we
obtain the following relation 
\begin{eqnarray}
B_{k,j;n,m}\left( x,y;q\right)
&=&[x]_{q}^{k}[y]_{q}^{j}\dsum\limits_{l=0}^{n}\dsum\limits_{r=0}^{m}\binom{n%
}{l}\binom{m}{r}  \notag \\
&&\times B_{l}^{\left( k\right) }\left( [1-x]_{q}\right) B_{r}^{\left(
j\right) }\left( [1-y]_{q}\right) \frac{\Delta ^{k}0^{n-l}}{k!}\frac{\Delta
^{j}0^{m-r}}{j!}  \label{Equation92}
\end{eqnarray}%
which is the relation of the $q$-Bernstein polynomials of two variables in
terms of Bernoulli polynomials of order $k$ and second Stirling numbers with
shift difference operator.

Let $\left( Eh\right) \left( x\right) =h\left( x+1\right) $ be the shift
operator. Then the $q$-difference operator is defined by 
\begin{equation}
\Delta _{q}^{n}=\dprod\limits_{j=0}^{n-1}\left( E-q^{j}I\right)
\label{Equation93}
\end{equation}%
where $I$ is and identity operator ( See \cite{Kim-Lee-Chae Jang} ).

For $f\in C[0,1]$ and $n\in 
\mathbb{N}
$, we have 
\begin{equation}
\Delta _{q}^{n}f\left( 0\right) =\dsum\limits_{k=0}^{n}\binom{n}{k}%
_{q}\left( -1\right) ^{k}q^{\binom{n}{2}}f\left( n-k\right) ,
\label{Equation94}
\end{equation}%
where$\binom{n}{k}_{q}$ is called the Gaussian binomial coefficients, which
are defined by 
\begin{equation}
\binom{n}{k}_{q}=\frac{[x]_{q}[x-1]_{q}\cdots \lbrack x-k+1]_{q}}{[k]_{q}!}.
\label{Equation95}
\end{equation}

\begin{theorem}
For $n,m,l,r\in 
\mathbb{N}
_{0}$, we have%
\begin{eqnarray*}
&&\frac{1}{\left( \left[ x\right] _{q}+\left[ 1-x\right] _{q}\right)
^{n-l}\left( \left[ y\right] _{q}+\left[ 1-y\right] _{q}\right) ^{m-l}}%
\dsum\limits_{k=l}^{n}\dsum\limits_{j=l}^{m}\frac{\binom{k}{l}\binom{j}{l}}{%
\binom{n}{l}\binom{m}{l}}B_{k,j;n,m}\left( x,y;q\right) \\
&=&\dsum\limits_{k=0}^{l}\dsum\limits_{j=0}^{l}q^{\binom{k}{2}+\binom{j}{2}}%
\binom{x}{k}\binom{y}{j}\left[ k\right] _{q}!\left[ j\right] _{q}!S\left(
l,k;q\right) S\left( l,j;q\right) .
\end{eqnarray*}

\begin{proof}
Let $F_{q}\left( t\right) $ be the generating function of the $q$-extension
of the second kind Stirling numbers as follows:%
\begin{equation}
F_{q}\left( t\right) :=\frac{q^{-\binom{k}{2}}}{[k]_{q}!}\dsum%
\limits_{j=0}^{k}\left( -1\right) ^{k-j}\binom{k}{j}_{q}q^{\binom{k-j}{2}%
}e^{[i]_{q}t}=\dsum\limits_{n=0}^{\infty }S\left( n,k;q\right) \frac{t^{n}}{%
n!}  \notag
\end{equation}%
From the above, we have 
\begin{equation}
S\left( n,k;q\right) =\frac{q^{-\binom{k}{2}}}{[k]_{q}!}\dsum%
\limits_{j=0}^{k}\left( -1\right) ^{j}q^{\binom{j}{2}}\binom{k}{j}%
_{q}[k-j]_{q}^{n}=\frac{q^{-\binom{k}{2}}}{[k]_{q}!}\Delta _{q}^{k}0^{n} 
\notag
\end{equation}%
where $[k]_{q}!=[k]_{q}[k-1]_{q}\cdots \lbrack 2]_{q}[1]_{q}.$ It is easy to
see that 
\begin{equation}
\lbrack x]_{q}^{n}=\dsum\limits_{k=0}^{n}q^{\binom{k}{2}}\binom{x}{k}%
_{q}[k]_{q}!S\left( n,k;q\right)  \label{Equation98}
\end{equation}%
by similar way 
\begin{equation}
\lbrack y]_{q}^{j}=\dsum\limits_{r=0}^{j}q^{\binom{r}{2}}\binom{y}{r}%
_{q}[r]_{q}!S\left( j,r;q\right) .  \label{Equation99}
\end{equation}%
We have above equality. Then, we obtain the desired result in Theorem from
the equations (\ref{Equation98}), (\ref{Equation99}) and Theorem 7.
\end{proof}
\end{theorem}

\section{Interpolation Function of Modified q-Bernstein Polynomial for
Functions of Two Variables}

For $s\in 
\mathbb{C}
$, and $x\neq 1$, $y\neq 1$, by applying the Mellin transformation to
generating function of Bernstein polynomials of two variables, we get%
\begin{eqnarray}
S_{q}\left( s,k,j;x,y\right) &=&\frac{1}{\Gamma \left( s\right) }%
\dint\limits_{0}^{\infty }F_{k,j}\left( -t,q;x,y\right) t^{s-k-j-1}dt  \notag
\\
&=&\frac{\left( -1\right) ^{k+j}[x]_{q}^{k}[y]_{q}^{j}}{k!j!}\left(
[1-x]_{q}+[1-y]_{q}\right) ^{-s}.  \label{Equation100}
\end{eqnarray}%
By using the equation (\ref{Equation100}), we define the interpolation
function of the polynomials $B_{k,j;n,m}\left( x,y;q\right) $ as follows:

\begin{definition}
Let $s\in 
\mathbb{C}
$ and $x\neq 1$, $y\neq 1,$ we define 
\begin{equation}
S_{q}\left( s,k,j;x,y\right) =\frac{[x]_{q}^{k}[y]_{q}^{j}}{k!j!}\left(
-1\right) ^{k+j}\left( [1-x]_{q}+[1-y]_{q}\right) ^{-s}.  \label{Equation50}
\end{equation}
\end{definition}

By using (\ref{Equation50}), we have $S_{q}\left( s,k,j;x,y\right)
\rightarrow S\left( s,k,j;x,y\right) $ as $q\rightarrow 1.$ Thus one has%
\begin{equation}
S\left( s,k,j;x,y\right) =\frac{\left( -1\right) ^{k+j}}{k!j!}%
x^{k}y^{j}\left( 2-\left( x+y\right) \right) ^{-s}.  \label{Equation51}
\end{equation}%
By substituting $x=1$ and $y=1$ into the above, we have $S\left(
s,k,j;x,y\right) =\infty $.

We now evaluate the $m$th $s$-derivatives of $S\left( s,k,j;x,y\right) $ as
follows:%
\begin{equation}
\frac{\partial ^{m}}{\partial s^{m}}S\left( s,k,j;x,y\right) =\log
^{m}\left( \frac{1}{2-\left( x+y\right) }\right) S\left( s,k,j;x,y\right)
\label{Equation52}
\end{equation}%
where $x\neq 1$ and $y\neq 1.$

\end{document}